\documentclass{amsart}
\usepackage {amssymb,latexsym}
\newtheorem{theorem}{Theorem}[section]

\newtheorem {proposition}[theorem]{Proposition}
\newtheorem {corollary} [theorem] {Corollary}

\theoremstyle{definition}

\theoremstyle{remark}
\newtheorem{remark}[theorem]{Remark}
\newtheorem{remarks}[theorem]{Remarks}
\begin{document}
\title{K-Theory Tools for Local and Asymptotic Cyclic Cohomology}

\author{Vahid Shirbisheh}
\address{Department of Mathematics and Computer Science, University of Tehran, Tehran, Iran}
\email{shirbish@khayam.ut.ac.ir \& shirbisheh@yahoo.com}

\date{February 27, 2002}

\begin{abstract}
A generalization of Connes-Thom isomorphism is given for stable, homotopy invariant, 
and split exact functors on separable $C^*$-algebras. As examples of these functors, 
we concentrate on asymptotic and local cyclic cohomology and the result is applied to improve 
some formulas in asymptotic and local cyclic cohomology of $C^*$-algebras. As an other application, it is shown that these cyclic theories are rigid after Rieffel's deformation quantizations.\\
\end{abstract}
\maketitle

\section*{Introduction}
Motivation for this paper comes from two sources. The first is Rieffel's and 
Abadie's papers [R6, A2], where it was proved that 
$K$-theory is rigid under Rieffel's deformation quantizations. The second is 
Rosenberg's works, [Ro3, Ro2]. In [Ro3], Rosenberg 
discussed the relations between $K$-theory and various quantization theories, 
and he suggested similar study on Connes' cyclic 
homology. In [Ro2], he studied behavior of algebraic $K$-theory under formal 
deformation quantization. Whereas, cyclic homology of formal deformation 
quantization has been studied by Nest and Tsygan in [NT1,NT2], it is natural to 
ask how cyclic (co)homology behaves under other 
quantization theories. Because of the algebraic nature of Connes' cyclic 
theory, we do not hope to get any satisfactory answer to this question on 
$C^*$-algebraic quantizations. Therefore, in such situations, we have to use 
other cyclic theories which do all duties of Connes' cyclic theory and are 
appropriate to deal with topological algebras. 

Bivariant local cyclic cohomology of admissible Fr\'echet algebras with 
compact supports and bivariant asymptotic cyclic 
cohomology of admissible Fr\'echet algebras are such theories, see [Pu1,Pu2,Pu3]. On 
commutative $C^*$-algebras, they are comparable with the cohomology of the character spaces as locally compact topological spaces. 
Moreover, they are stable and homotopy invariant 
bifunctors which appropriate excisions are hold in their both variables and  there 
are bivariant Chern-Connes characters from $KK$-theory into them.

Among $C^*$-algebraic quantization theories, Rieffel's deformation 
quantizations, [R3, R5], are the best theories for our 
purpose. Because they are associated easily to crossed product algebras (up to 
strong Morita equivalence), [R6, A2], and the set 
of their examples contains several important noncommutative spaces.

It is clear from [R6,A2] that in order to repeat their proofs for other 
functors, we need Morita invariance, Connes-Thom 
isomorphism, Pimsner-Voiculescu exact sequence, and Bott periodicity for them. 
So, we provide  these tools in section 1. 
Section 2 is devoted to the application of these tools in local and specially 
asymptotic cyclic cohomology. In section 3 we show 
local and asymptotic cyclic (co)homology groups of noncommutative and 
commutative Heisenberg manifolds are isomorphic. As another example of strict 
deformation quantization, we study behavior of local and asymptotic cyclic 
theories under deformation quantization by action of $\mathbb{R}^n$.
\par
In this paper, we restrict ourselves to separable $C^*$-algebras and their 
dense $^*$-subalgebras. $\mathcal{K}$ denotes the $C^*$-algebra of compact 
operators on an infinite dimensional separable Hilbert space. The functor $S(-)$ 
is suspension, i.e. for any $C^*$-algebra $A, S(A)=A\otimes C_0(\mathbb{R})$. A 
covariant functor is said simply {\it functor} and what will be 
proved for functors are hold for contravariant functors, {\it cofunctors}, by 
similar proofs. Let {\it F} be a functor, it is called {\it stable} if for 
every $C^*$-algebra $A$ the natural embedding $A\rightarrow A\otimes\mathcal{K}$ 
induces a natural isomorphism $F(A)\rightarrow F(A\otimes\mathcal{K})$, $F$ is 
called {\it homotopy invariant} if $F(f)=F(g):F(A)\rightarrow F(B)$ whenever 
$f,g\in Hom (A,B)$ are homotopic and $F$ is called {\it split exact} when for 
every split exact sequence 
$0 \rightarrow A \stackrel {f}{\rightarrow} B \stackrel {p}{\rightarrow} C \rightarrow 0$
with splitting *-homomorphism $g:C\rightarrow B$, the map 
$F(f)\oplus F(g):F(A)\oplus F(C)\rightarrow F(B)$ is an isomorphism. 

By "{\it bivariant local cyclic cohomology}" ("{\it bivariant asymptotic cyclic 
cohomology"}), we mean bivariant local cyclic cohomology of admissible Fr\'echet 
algebras with compact supports (bivariant asymptotic cyclic cohomology of 
admissible Fr\'echet algebras), and for a pair of $C^*$-algebras (A,B) we denote 
it by $HC_{lc}^*(A,B)$ ($HC_{\alpha}^*(A,B)$). We do similarly for local 
(asymptotic) cyclic homology and cohomology groups. 
{\it Bivariant Chern-Connes character} from $KK$-theory to bivariant local and 
asymptotic cyclic cohomology, as defined in [Pu1, Pu3], are denoted 
respectively as follows
\[
ch_{biv}:KK^*(-,-)\rightarrow HC_{lc}^*(-,-)
\]
and 
\[
ch_{biv}^s:KK^*(-,-)\rightarrow HC_{\alpha}^*(S-,S-).
\]

\section{connes-thom isomorphism}
We use some techniques of $KK$-theory to prove Connes-Thom isomorphism. Our 
approach is known as Cuntz's picture of $KK$-theory, [C3]. We note Bott 
periodicity and Pimsner-Voiculescu exact sequence were 
previously studied by the same way in [C2] and our 
proof for theorem 1.2 is an application of Fack and Skandalis result, [FS].  

A {\it quasihomomorphism} between two $C^*$-algebras $A$ and $B$ is a diagram 
as follows:
\[A \stackrel {\alpha}{\underset{\bar{\alpha}}{\rightrightarrows}} E\rhd J\stackrel{\mu}{\rightarrow}B,
\]
where $E$ and $J$ are $C^*$-algebras, $J \lhd E$, and $\alpha, \bar\alpha,\mu$ are 
$^*$-homomorphisms such that
\begin{itemize}
\item[(i)] 
$\mu$ is an inclusion,
\item[(ii)] 
$E$ is the $C^*$-algebras generated by $\alpha(A)$ and $\bar{\alpha}(A)$,
\item[(iii)] 
$J$ is the closed two-sided ideal generated by $\alpha(x)-\bar{\alpha}(x), x\in A$ in 
$E$,
\item[(iv)] the composition of $\alpha$ and the quotient map 
$E\rightarrow {E}/{J}$ is injective, thus, an isomorphism.
\end{itemize}
(A diagram as above with only the property 
$\alpha(x)-\bar{\alpha} (x)\in J,\text {for}\; x\in A$ is 
called {\it prequasihomomorphism}.) Now, let $A$ be a $C^*$-algebras, by 
definition, $QA$ is the universal $^*$-algebra generated by symbols 
$x, q(x), x\in A$ satisfying in the relation 
$q(xy)=xq(y)+q(x)y-q(x)q(y)$ and let $qA$ is defined as the ideal generated by 
$q(x), x\in A$ in $QA$. We equip $QA$ with the largest $C^*$-norm, 
\[
\|x\|_{\infty}=\sup\{\|\pi(x)\|; \pi\;\;\text {is a}\ ^*-\text {representation}\}
\]
and still write $QA$, $qA$ for completions of $QA$, $qA$ with respect to this 
norm. 

For two $C^*$-algebras $A$ and $B$, let $[qA, B \otimes \mathcal{K}]$ be the 
set of homotopy classes of $^*$-homomorphisms from $qA$ to 
$B \otimes \mathcal{K}$, on this set, addition is defined by 
$[\varphi]+[\psi]\!=\!\left[\begin{pmatrix}
\varphi&0\\
0& \psi
\end{pmatrix}\right]$. 
By this addition it becomes an Abelian group equal to $KK(A,B)$, (or more 
exactly, $KK_0(A,B)$), for details see [C1, C3].

We have two homomorphisms $\imath,\;\bar{\imath}$ from $A$ into $QA$ defined by 
$\imath(x)=x$, $\bar{\imath}=x-q(x)$. The quasihomomorphism given by the 
diagram 
\[A \underset {\bar {\imath}} {\stackrel {\imath} {\rightrightarrows}} QA \rhd qA
\]
is the {\it universal quasihomomorphism} from $A$ into $qA$. Let 
$\varphi\!:\!qA\!\rightarrow\!\!B$ be a $^*$-homomorphism, it is naturally 
associated to the quasihomomorphism
\[
A\overset{\pi\imath}{\underset{\pi\bar{\imath}}{\rightrightarrows}}{QA}/{\ker\varphi}\rhd
{qA}/{\ker\varphi}\;\;\overset{\bar\varphi}{\rightarrow} B,
\]
where $\pi:QA\rightarrow 
{QA}/{\ker\varphi},\;\;\bar\varphi:{qA}/{\ker\varphi}\rightarrow B$ 
are respectively, the quotient map, and the inclusion 
map defined by $\varphi$. (and vice versa, every quasihomomorphism from $A$ to $B$ give rise to a 
homomorphism from $qA$ to $B$, for details see proposition 1.1 of [C3].) 

Following result is the main tool in our study. It was stated for functors 
from the category of $C^*$-algebras to the category 
of $\mathbb{Z}$-modules
 in [C3], and can be restated for any commutative ring ${R}$, 
instated of $\mathbb{Z}$. 

\begin{proposition} (Cuntz, [C3, 2.2.a]) {\it Let $F$ be a stable, 
homotopy invariant, split exact functor from the category of 
$C^*$-algebras to the one of $R$-modules, every} $\varphi\in KK(A, B)$ 
{\it induces a morphism} $F(\varphi):F(A)\rightarrow F(B)$ {\it compatible
 with the Kasparov product, i.e.} $F(\psi\varphi)=F(\psi)F(\varphi)$ for any $\psi\in KK(B, C)$.
\end {proposition}
\begin{proof} We sketch only the definition of $F(\varphi)$. Consider 
$\varphi$ as a $^*$-homomorphism from $qA$ to $B$, as we saw this 
$^*$-homomorphism is associated to the quasihomomorphism 
\[A \underset {\pi \bar {\imath}} {\stackrel {\pi\imath} {\rightrightarrows}} {QA}/\ker\varphi\rhd{qA}/\ker\varphi \stackrel {\bar\varphi}{\rightarrow}B,
\]
and this quasihomomorphism induces the desired morphism as 
$F(\varphi)=F(\bar{\varphi})\big(F(\pi\imath)-F(\pi\bar{\imath})\big)$. 
According to the definition, it is obvious that 
$F(\varphi)(r-)=r\ F(\varphi)(-)$ for $r\in R$. 
\end{proof}
Above result shows that in order to obtain a desired morphism between $F(A)$ 
and $F(B)$ it is enough to find an appropriate element of $KK(A,B)$.
 An element of $KK(A,B)$ inducing isomorphism is called 
{\it $KK$-equivalence}. Equivalently, $\bold{x}\in KK(A, B)$ is a 
$KK$-equivalence, if there is $\bold{y}\in KK(B , A)$ such that 
$\bold{xy}=1_B,\bold{yx}=1_A$. If there is a $KK$-equivalence in $KK(A,B)$, $A$ 
and $B$ are called {\it $KK$-equivalent}. \par
In [C2] Bott periodicity and Pimsner-voiculescu exact sequence were made for functors described in 1.1. Now 
Connes-Thom isomorphism is accessible by applying $KK$-equivalence 
$\bold{t}_{\alpha}\in KK(A\times_{\alpha}\mathbb{R}, SA)$ made in section 19.3 of [B].
\begin {theorem} {Let $F$ be a functor as 1.1, and 
$\alpha$ be a (strongly) continuous action of $\Bbb{R}$ on a $C^*$-algebra $A$, 
then $F(\bold{t}_{\alpha})$ is an isomorphism between 
$F(A\times_{\alpha}\Bbb{R})$ and $F(SA)$}.
\end {theorem}
\begin {remarks}
\begin{itemize}
\item[(a)]
{Indeed, theorem 1.2 is not exactly the generalization of Conn- es' isomorphism. However, if we define $F_{-i}(A)=F(S^iA),\ i\in\Bbb{Z}$, then we have Connes' isomorphism, $F_*(A\times_{\alpha}\Bbb{R})\cong F_{*-1}(A)$. In $K$-theory, $K$-homology and other general $K$-functors like $KK(B_1,-\otimes B_2)$, the assumption $F_{-i}=F(S^iA)$ is automatically hold, but, we still do not know whether it is true about local and asymptotic cyclic theories. So, we have to content ourselves with the theorem 1.2 at present, see next section for full Connes-Thom isomorphism for asymptotic and local cyclic cohomology.} 
\item[(b)] 
$KK$-equivalence elements of $KK(A,B)$ preserve both torsion and torsion-free parts of $K$-groups, while if the ring $R$ contains $\Bbb{Q}$, one can consider invertible elements (with respect to the Kasparov product) of the module $KK(A, B)\otimes_{\Bbb{Z}}R$ which preserve only torsion-free parts of $K$-groups, thus, we have other choices to construct new equivalences.
\item[(c)] 
$K$-amenability; locally compact group $G$ is called $K$-amenable if $\bar{p}$, the element of $KK\!\big(\!C^*(\!G),\!C^*_r(\!G)\!\big)$ induced by the projection $p:C^*(G)\rightarrow C^*_r(G)$, is a $KK$-equivalence. So $K$-theory of group $C^*$-algebra of a $K$-amenable group equals $K$-theory of its reduced group $C^*$-algebra. Similarly, we can define $K$-amenability for any arbitrary functor $F$ Also, for any commutative ring $R$, we can define $K(R)$-amenability, i.e. $G$ is $K(R)$-amenable, whenever $\bar{p}\in KK\big(C^*(G),C_r^*(G)\big)\otimes_{\Bbb{Z}}R$ be an isomorphism.
\item[(d)] 
Let $\alpha$ be an automorphism of a $C^*$-algebra $A$. As it was done in [C2], 
Pimsner-Voiculescu exact sequence can be constructed using a 
$KK\!$ -equival- ence in $KK(A,T_{\alpha})$, where $T_{\alpha}$ is Toeplitz 
algebra associated to the automorphism $\alpha$. This $KK\!$-equivalence
 allows us to replace $K$-groups of Toeplitz algebra with $K$-groups 
of $A$ in the six term exact sequence associated to the following short exact 
sequence, known as Toeplitz extension: 
\[
0\rightarrow \mathcal{K}\otimes A\rightarrow T_{\alpha}\rightarrow 
A\times_{\alpha}\Bbb{Z}\rightarrow 0 
\]
Therefore, in order to obtain P-V exact sequence for a functor, only conditions 
assumed in 1.1 is required. \par
Since Pimsner's and Vioculescu's paper, [PV], was appeared, their work has found several generalizations, [S,E2,KhS,P,AEE]. In 
all of them $A\times_{\alpha}\Bbb{Z}$ is replaced by a new $C^*$-algebra, 
e.g. $A\bar{\times}_E\Bbb{Z}$, which generalizes the crossed product of 
$C^*$-algebra $A$ by $\Bbb{Z}$ and we have a generalized Toeplitz extension of 
$A$ by a generalized Toeplitz algebra, e.g. $\bar{T}_E$, which is 
$KK$-equivalent to $A$. The six term exact sequences of such extensions are 
considered as generalizations of P-V exact sequence. One can repeat 
our discussion for these generalizations. 
\end{itemize}
\end {remarks}

Let $\alpha$ be an automorphism of a $C^*$-algebra $A$, the mapping torus of 
$\alpha$ is 
\[
M_{\alpha}=\{f:[0,1]\rightarrow A; f(1)=\alpha(f(0))\}.
\]
Consider the action of $\Bbb{Z}$ on $A$ defined by $\alpha$, its dual $\hat{\alpha}$ 
is an action of the dual of $\Bbb{Z}$, $\Bbb{T}=\hat{\Bbb{Z}}$, on 
$A\times_{\alpha}\Bbb{Z}$. Let $\pi:\Bbb{R}\rightarrow{\Bbb{R}}/{\Bbb{Z}}$ 
be the quotient map, then $\alpha'=\hat{\alpha}o\pi$ is an action of $\Bbb{R}$ on 
$A\times_{\alpha}\Bbb{Z}$, which is trivial on $\Bbb{Z}$ and 
$A\times_{\alpha}\Bbb{Z}\times_{\alpha'}\Bbb{R}$
 is isomorphic to the mapping torus of $\hat{\hat{\alpha}}$ on 
$A\times_{\alpha}\Bbb{Z}\times_{\hat{\alpha}}\Bbb{T}$,
 (see [B], proposition 10.3.2). 
By Takai duality, it means $A\times_{\alpha}\Bbb{Z}\times_{\alpha'}\Bbb{R}$
is isomorphic to the mapping torus of $\alpha\otimes Ad_{\rho}$ on $A\otimes\mathcal{K}$, 
where $\rho$ is the right regular representation of $\Bbb{Z}$ and $\mathcal{K}$ is 
thought of as the $C^*$-algebra of 
compact operators on $\ell^2(\Bbb{Z})$. Now, suppose $\alpha,\beta$ be two homotopic 
automorphism of $A$, then $\alpha\otimes Ad_{\rho}$ and $\beta\otimes Ad_{\rho}$ are 
homotopic too. Applying proposition 10.5.1 of [B] we deduce 
$A\times_{\alpha}\Bbb{Z}\times_{\alpha'}\Bbb{R}$ is isomorphic to 
$A\times_{\beta}\Bbb{Z}\times_{\beta'}\Bbb{R}$. Thus, as a generalization of corollary 10.5.2 
of [B], we have the following result:
\begin{corollary} Let $F$ be a functor as 1.1 and $\alpha,\beta$ be two homotopic automorphism of a $C^*$-algebra $A$, then 
$F(A\times_{\alpha}\Bbb{Z})\cong F(A\times_{\beta}\Bbb{Z})$.
\end {corollary}
\subsection* {Morita invariance}
Algebraic Morita equivalence provides the opportunity to replace an algebra 
$A$ with another algebra $B$ to simplify computations 
of (co)homology groups. For example, cyclic and Hochschild (co)homology and 
algebraic $K$-theory have Morita invariance 
property, [L, Ro1]. Since strong Morita equivalence for operator algebras was 
defined by Rieffel in [R1, R2], it has found 
several examples and applications. Now, it is an usual approximation for 
noncommutative spaces. For instance , $K$-theory of two strongly Morita 
equivalent $C^*$-algebras are isomorphic, [E1]. It is well known 
that two strongly Morita equivalent $\sigma$-unital $C^*$-algebras are 
stably isomorphic, [BGR]. Since every separable $C^*$-algebra is $\sigma$-unital and because of our assumptions, we have
\begin {remark}
In our discussion, all functors have strong Morita invariance 
property on separable $C^*$-algebras.
\end {remark}

\section{Local and Asymptotic Cyclic Cohomology}
Some applications of $K$-theory of operator algebras are related to the 
structural problems on $C^*$-algebras, e.g. classification theories, while 
others are generalizations of well known problems on topological $K$-theory for 
noncommutative spaces, for instance, index theorems. Since $K$-theory is a 
powerful functor with a number of theorems and techniques, and also it is defined 
for all $C^*$- and Von Neumann algebras and their spectral invariant dense 
subalgebras, it seems other functors can not do any more about the first part 
of applications. In the second part of applications, there were some 
questions on existence of stable homology and cohomology theories on noncommutative 
spaces, $C^*$-algebras, similar to homology and cohomology theories on locally 
compact topological spaces, as commutative spaces. Appearance of asymptotic and 
local cyclic cohomology was an answer for these questions. As we have developed 
some tools of $K$-theory for a class of functors containing local and 
asymptotic cyclic theories, it is the time we study some of their applications.

Asymptotic and local cyclic cohomology are stable in both variables, (see 
theorem 8.18 of [Pu1] and corollary 4.10 of [Pu2]). Also, composition products 
and continuous homotopy theorems prove their homotopy invariances, (see theorem 
6.5 and 6.15 of [Pu1] and theorem 3.5 and 3.18 of [Pu2]). But, there is a 
difference between local cyclic theory and asymptotic cyclic theory about split 
exactness. First we study asymptotic cyclic cohomology.\\
\subsection*{Asymptotic cyclic cohomology} At the moment, we have only stable 
excisions for asymptotic cyclic homology and cohomology, i.e. for any short exact sequence 
of separable $C^*$-algebras
\[
0\rightarrow I \overset{i}{\rightarrow} A\overset{s}{\underset{f}{\leftrightarrows}} B\rightarrow 0
\]
with a bounded linear section $s:B \rightarrow A$, and any $C^*$-algebra $C$, there are 
six term exact sequences 
\begin{equation*}\begin{matrix}
HC_{\alpha}^0(C,SI) & \longrightarrow & HC_{\alpha}^0(C,SA) & \overset{Sf_*}
{\longrightarrow} & HC_{\alpha}^0(C,SB)\\
\uparrow\partial & & & & \downarrow\partial \\
HC_{\alpha}^1(C,SB) & \overset{Sf_*}{\longleftarrow} & HC_{\alpha}^1(C,SA) & 
\longleftarrow & HC_{\alpha}^1(C,SI)
\end{matrix}\ 
\end{equation*}
and
\begin{equation*}\begin{matrix}
HC_{\alpha}^0(SI,C) & \longleftarrow & HC_{\alpha}^0(SA,C) & \overset{Sf^*}
{\longleftarrow} & HC_{\alpha}^0(SB,C)\\
\downarrow\partial & & & & \uparrow\partial \\
HC_{\alpha}^1(SB,C) & \overset{Sf^*}{\longrightarrow} & HC_{\alpha}^1(SA,C) & 
\longrightarrow & HC_{\alpha}^1(SI,C)
\end{matrix}\ 
\end{equation*}
\\
If the given short exact sequence of $C^*$-algebras is split, these six term 
exact sequences give rise to four split short exact sequences of complex vector 
spaces as follows:
\[
0\rightarrow HC_{\alpha}^*(C, SI)\rightarrow HC_{\alpha}^*(C, SA)\rightarrow HC_{\alpha}^*(C, SB)\rightarrow 0, *=0, 1,
\]
\[
0\rightarrow HC_{\alpha}^*(SB, C)\rightarrow HC_{\alpha}^*(SA, C)\rightarrow HC_{\alpha}^*(SI, C)\rightarrow 0, *=0, 1.
\]
\\
Therefore, for any $C^*$-algebra $C$, functor $HC_{\alpha}^*(C, S-)$ and cofunctor 
$HC_{\alpha}^*(S-, C)$ are split exact, thus, we have the following results:
\\
\begin {theorem} Let $A, B$ be two separable $C^*$-algebras and 
$\alpha$ be a (strongly) continuous action of $\Bbb{R}$ on $A$, then
\\
\begin{enumerate}
\item $HC_{\alpha}^*(A, B)\cong HC_{\alpha}^{*+1}(SA, B)\cong HC_{\alpha}^{*+1}(A, SB)$,
\\
\item $HC_{\alpha}^*(A, B)\cong HC_{\alpha}^*(SA, SB)$,
\\
\item $HC_{\alpha}^*(A\times_{\alpha}\Bbb{R}, B)\cong HC_{\alpha}^{*+1}(A, B)$,
\\
\item $HC_{\alpha}^*(B, A\times_{\alpha}\Bbb{R})\cong HC_{\alpha}^{*+1}(B, A).$
\\
\end{enumerate}
\end {theorem}
\begin {proof} For a given $C^*$-algebra $D$, suppose    
$\alpha_{SD}\in HC_{\alpha}^1(S^2D, SD)$ and $\beta_{SD}\in HC_{\alpha}^1(SD, S^2D)$ 
be respectively Dirac and Bott elements defined in definition 9.3 of [Pu1]. 
Also let $\hat{i}$ be the dual action of the trivial action of $\Bbb{R}$ on 
$D$, then
\\
\begin{align*}
HC_{\alpha}^*(A, B)&\cong HC_{\alpha}^*\big(S(A\times_{\hat{i}}\Bbb{R}), B\big)\\
&\cong HC_{\alpha}^*(S^2A, B)\\
&\cong HC_{\alpha}^{*+1}(SA, B),
\end{align*}
where isomorphisms respectively come from Takai duality, theorem 1.2 and 
composition by $\beta_{SA}$. Similarly, the second isomorphism of (1) is proved. From (1), 
(2) is obvious. For (3) we observe
\begin{align*}
HC_{\alpha}^*(A\times_{\alpha}\Bbb{R},
 B)&\cong HC_{\alpha}^{*+1}\big(S(A\times_{\alpha}\Bbb{R}),B
\big)\\
&\cong HC_{\alpha}^{*+1}(S^2A, B)\\
&\cong HC_{\alpha}^{*+1}(A, B),
\end{align*}
and similarly for (4).\end {proof}
     
Parts (3) and (4) are full Connes-Thom isomorphism for asymptotic cyclic 
cohomology and homology, which we promised in previous section. Using above 
theorem we can state unstable (third!) excision theorem for asymptotic cyclic 
cohomology, see also [Pu1]. 
\\ 
\begin {theorem} For any $C^*$-algebra $C$ and any short exact 
sequence of $C^*$-algebras
\[
0\rightarrow I\overset{i}{\rightarrow} A\overset{s}{\underset{f}{\leftrightarrows}} B\rightarrow 0
\]
with a bounded linear section $s: B\rightarrow A$, there exist following six term exact 
sequences:
\begin{equation*}\begin{matrix}
HC_{\alpha}^0(C,I) & \longrightarrow & HC_{\alpha}^0(C,A) & \overset{\hat{f}_*}
{\longrightarrow} & HC_{\alpha}^0(C, B)\\
\uparrow\partial '& & & & \downarrow\partial '\\
HC_{\alpha}^1(C,I) & \overset{\hat{f}_*}{\longleftarrow} & HC_{\alpha}^1(C,A) & 
\longleftarrow & HC_{\alpha}^1(C, B)
\end{matrix}\
\end{equation*}
\\
\begin{equation*}\begin{matrix}
HC_{\alpha}^0(I,C) & \longleftarrow & HC_{\alpha}^0(A,C) & \overset{\hat{f}^*}
{\longleftarrow} & HC_{\alpha}^0(B,C)\\
\downarrow\partial '& & & & \uparrow\partial '\\
HC_{\alpha}^1(B,C) & \overset{\hat{f}^*}{\longrightarrow} & HC_{\alpha}^1(A,C) & 
\longrightarrow & HC_{\alpha}^1(I,C)
\end{matrix}\
\end{equation*}
\end {theorem}
\noindent{\bf Notes:} Maps $\partial ', \hat{f}_*, \hat{f}^*$ are appropriate compositions 
of maps $\partial , Sf_*, Sf^*$ and isomorphisms employed in the proof of theorem 
2.1. Also, one can consider another Chern character form $KK$-theory into 
unstable bivariant asymptotic cyclic cohomology. Of course, we do not know 
whether the new Chern character is natural.
\par
As another consequence of theorem 2.1 we have the following corollary:
\\
\begin {corollary} 
\begin {enumerate} \item Let $X$, $Y$ be two finite $CW$-complexes, 
then 
\begin{equation*}
HC_{\alpha}^*\big(C(X), C(Y)\big)\cong\underset{\pmod{2}}{\underset{n+m\equiv *}
{\text {Hom}}}\big(\overset{\infty}{\underset{n=0}{\oplus}}H^n(X,\Bbb{C}), 
\overset{\infty}{\underset{m=0}{\oplus}}H^m(Y,\Bbb{C})\big), 
\end{equation*}
where $H^*(X,\Bbb{C})$ denotes singular cohomology of $X$ with coefficients in 
$\Bbb{C}$.\\ 
\item Let $X$ be a locally compact metrisable topological space (or equivalently, 
$C(X)$ is a separable $C^*$-algebra), then
\[
HC^{\alpha}_*\big(C(X)\big)=HC_{\alpha}^*\big(\Bbb{C}, C(X)\big)\cong\overset{\infty}
{\underset{n=-\infty}{\oplus}}H_c^{*+2n}(X,\Bbb{C}),
\]
where $H_c^*(X,\Bbb{C})$ denotes sheaf cohomology with compact supports and 
coefficients in $\Bbb{C}$.\\
\end {enumerate}
\end {corollary}
\begin{proof} See theorems 11.2 and 11.7 of [Pu1].
\end {proof}

\subsection*{Local cyclic cohomology}
All of aforementioned results on asymptotic cyclic cohomology can be proved (even easier) for local 
cyclic cohomology too. Therefore, we bring only some remarks here.

\begin {remarks} 
\begin{itemize}
\item[(a)] 
Since we have a natural transformation compatible 
with composition products from bivariant asymptotic cyclic 
cohomology to bivariant local cyclic cohomology, all asymptotic cyclic 
equivalences like $\alpha_{S-},and \beta_{S-}$, 
induce similar local cyclic equivalences, see 11.9.b of [Pu1]
 and 3.23 of [Pu2]. 
\item[(b)] 
Already there is an unstable excision for local cyclic (co)homology, see 
theorem 5.12 of [Pu3]. 
\end{itemize}
\end {remarks}

Since local cyclic cohomology behave reasonably under inductive limits, corollary 
2.3 is adapted as follows for local cyclic theory:\\
\begin {corollary} Let $A$, $B$ be separable, commutative 
$C^*$-algebras with corresponding locally compact spaces $X$, $Y$, then 
\begin{equation*}
HC_{lc}^*(A,B)\cong\underset{\pmod{2}}{\underset{n+m\equiv *}
{\text {Hom}}}\big(\overset{\infty}{\underset{n=0}{\oplus}}H^n_c(X,\Bbb{C}), 
\overset{\infty}{\underset{m=0}{\oplus}}H^m_c(Y,\Bbb{C})\big)\ .
\end{equation*}
\end {corollary}

\begin {remarks}
\begin{itemize}
\item[(a)] 
Local cyclic homology groups of commutative separable $C^*$-algebras are 
isomorphic to their local cyclic cohomology groups. This is true more generally 
on the class of separable $C^*$-algebras which are strong Morita equivalent to 
commutative separable $C^*$-algebras. 
\item[(b)] 
let $\mathcal{C}$ be the class of $C^*$-algebras described in the theorem 10.7 
of [Pu1], then $ch_{biv}$,and $ch^s_{biv}$ yield following isomorphisms for any $A$, 
$B$ in $\mathcal{C}$:
\\
\begin{align*}
&ch: KK^*(A,B)\otimes_{\Bbb{Z}}\Bbb{C}\overset{\cong}{\longrightarrow}HC^*_
{lc}(A,B), 
\\
&ch: KK^*(A,B)\otimes_{\Bbb{Z}}\Bbb{C}\overset{\cong}{\longrightarrow}HC^*_
{\alpha}(A,B).
\end{align*}
\\
By proposition 6.2 of [C2] all commutative separable $C^*$-algebras belong to 
$\mathcal{C}$, thus, torsion free part of $K$-theory groups and $K$-homology
 groups of commutative $C^*$-algebras are isomorphic.
\item[(c)] 
As an example we have: 
$HC^{\alpha}_*\big(C(\Bbb{T}^n)\big)\cong HC_{\alpha}^*\big(C(\Bbb{T}^n)\big)
\cong\Bbb{C}^{2^{n-1}},\; *=0,1$, 
$HC^{lc}_*\big(C(\Bbb{T}^n)\big)\cong HC_{lc}^*\big(C(\Bbb{T}^n)\big)
\cong\Bbb{C}^{2^{n-1}},\; *=0,1$, see [R4].
\end{itemize}
\end {remarks}
\section {Strict Deformation Quantization}
Motivated by formal deformation quantization, Rieffel in [R3] introduced a $C^*$-algebraic 
framework for deformation quantization known as strict deformation quantization. Suppose $A$ be a $C^*$-algebra with a dense $^*$-subalgebra $\mathcal{A}$ equipped with a Poisson bracket $\{,\}$, {\it a strict deformation quantization} of $A$ in the direction of $\{,\}$ consists of an open Interval $I$ containing $0$, and a family of pre-$C^*$-algebra structures $\{(\times_{\hbar}, ^*\!\,_{\hbar}, ||\;\;\|_{\hbar})\}_{\hbar\in I}$  on $\mathcal{A}$ which for $\hbar=0$ it coincide with the pre-$C^*$-structure inherited from $A$ such that if we denote by $A_\hbar$ the completion of $\mathcal{A}$ under $C^*$-structure $(\times_\hbar, ^*\!\,_\hbar, \|\;\;\|_\hbar)$, then the family $\{A_\hbar\}_{\hbar\in I}$ with constant sections from $I$ into $\mathcal{A}$ constitutes a continuous field of $C^*$-algebras and for every $a,b\in\mathcal{A}$ we have 
\[
\underset{\hbar\rightarrow 0}{\lim}\|\frac{(a\times_{\hbar}b-ab)}{i\hbar}-\{a,b\}\|=0.
\]
The first examples of this definition are noncommutative tori and noncommutative Heisenberg manifolds. Deformation quantization of $C^*$-algebras by actions of finite dimensional real vector spaces, formulated in [R5], provides other examples for strict deformation quantization. In these constructions, the actions of $\Bbb{Z}$ and $\Bbb{R}$ play important roles As a consequence, the quantized algebras naturally are associated to appropriate crossed product $C^*$-algebras by actions of $\Bbb{Z}$ and $\Bbb{R}$. Thus, Connes-Thom isomorphism and Pimsner-Voiculescu exact sequence can be used to compare $K$-theory of the quantized $C^*$-algebra with the one of the original $C^*$-algebra, [R6,A2]. Now, these tools have been generalized for stable, split exact, and homotopy invariant functors. So, we can repeat Rieffel's and Abadie's works to obtain the same results for these functors, in particular for local and asymptotic cyclic (co)homology.
\\
\subsection {Noncommutative Heisenberg Manifolds}
For each positive integer $c$, the Heisenberg manifold $M_c$ is defined by the quotient ${G}/{D_c}$, where $G$ is the 
Heisenberg group, 
\[
G=\left\{\begin {pmatrix}
1&y&z\\
0&1&x\\
0&0&1
\end {pmatrix}; x,y, z \in \Bbb{R}\right\}\,
\]
\\
and
$D_c$ is the subgroup of G with $x , y, cz\in\Bbb{Z}$. Every non-zero Poisson bracket on $M_c$ is determined by two real numbers $\mu, \nu$, where $\mu^2+\nu^2\neq0$. According to [A2], deformation quantization of $M_c$ in the direction of non-zero Poisson bracket $(\mu, \nu)$, which is invariant under the action of $G$ by left translation, is denoted by 
$\{D^{c,\hbar}_{\mu\nu}\}_{\hbar\in\Bbb{R}}$, where $\hbar$ is the parameter of deformation. We know from [R3,A2] that for $\hbar\neq0$, the algebra 
$D^{c,\hbar}_{\mu\nu}$ is the generalized fixed-point algebra of 
$C_0(\Bbb{R}\times\Bbb{T})\times_{\lambda^\hbar}\Bbb{Z}$ under the action $\rho$, where
\[
\lambda^{\hbar}_k(x,y)=(x+2k\hbar\mu, y+2k\hbar\nu),\quad k\in\Bbb{Z},
\]
and if $e(x)=\exp(2\pi ix)$, then
\[
(\rho_k\Phi)(x,y,p)=e\big(ckp(y-{\hbar}p\nu)\big)\Phi(x+k, y,p), k\in\Bbb{Z}.
\]
It follows from theorem 2.11 of [A2] that the algebra $D_{\mu\nu}^{c,\hbar}$ is strongly Morita equivalent to another fixed-point algebra obtained by the action $\gamma^{\lambda^{\hbar}}$ of $\Bbb{Z}$ on 
$C_0(\Bbb{R}\times\Bbb{T})\times_{\sigma}\Bbb{Z}$, where
\[
\sigma_k(x,y)=(x-k,y),\quad k\in\Bbb{Z},
\]
and 
\[
(\gamma_p^{\lambda^{\hbar}}\Phi)(x,y,k)=e\big(-ckp(y-{\hbar}p\nu)\big)\Phi(x-2p{\hbar}\mu,
y-2p{\hbar}\nu,k),\;\;p\in\Bbb{Z}.
\]
Now, corollary 1.4 and homotopy ${\hbar}\rightarrow\lambda^{\hbar}$ show local and asymptotic cyclic (co)homology groups of $D^{c,\hbar}_{\mu\nu}$ are independed of parameter $\hbar$. On the other hand for any real number $\hbar, D^{c,\hbar}_{\mu\nu}$ and 
$D^{c,1}_{\hbar\mu,\hbar\nu}$ are isomorphic. So, one can drop parameter $\hbar$ from the notation and simply write $D^c_{\mu,\nu}$ instead of $D^{c,1}_{\hbar\mu,\hbar\nu}$. 
Also, for any pair of integers $k,l$, $D^{c,\hbar}_{\mu\nu}$ and 
$D^{c,\hbar}_{\mu+k,\nu+l}$ are isomorphic, see proposition 1 of [A1]. Thus, the assumption $\hbar\neq0$ can be ignored. This shows every stable, homotopy invariant, and split exact (co)functor $F$ is rigid under deformation quantization of Heisenberg manifolds. As a consequence, we have
\begin {theorem} For any real numbers $\mu,\nu,\hbar$ and integer $c$, even and odd asymptotic and local cyclic homology and cohomology groups of Heisenberg manifold $D^{c,\hbar}_{\mu\nu}$ are isomorphic to $\Bbb{C}^3$.
\end {theorem}
\begin {proof}
It is enough to consider only the commutative, which is an easy application of remark 2.6(b) and theorem 3.4 of [A2].
\end {proof}
\subsection{Deformation Quantization of $C^*$-Algebras by Actions of $\mathbb {R}^n$}
Let $\alpha$ be a strongly continuous action of $\Bbb{R}^n$ on $C^*$-algebra
 $A$, and $A^{\infty}$ be the dense $^*$-subalgebra of its smooth vectors. Also, let $J$ be a skew-symmetric matrix on $\Bbb{R}^n$. \\
On $A^{\infty}$ deformed product $\times_{J}$ is defined by
\[
a\times_Jb=\iint\alpha_{Ju}(a)\alpha_v(b)e(u\cdot v)\quad a,b\in A^{\infty},
\]
where $e(t)=\exp(2\pi it)$. Also, a $C^*$-norm $\|\;\;\|_J$ and an involution 
$*_J$ compatible with the product $\times_J$ are defined on $A^{\infty}$,(for details see [R5]). 
The completion of this pre-$C^*$-structure is denoted by $A_J$ and called \emph {quantization of $A$ by the action of $\alpha$ in the direction of $J$. 
Let $A$} be separable, then it is $\sigma$-unital, and consequently 
$A_J$ is $\sigma$-unital too, see [R5]. It was shown in [R6] that $A_J$ is strong Morita equivalent to a stable crossed product 
of $A$ as follows:
\[
A_J\overset{M}{\simeq}A\times_{\rho}\Bbb{R}^n\otimes
 C_0(\Bbb{R}^m)\otimes \mathcal{K},
\]
where $m$ is the dimension of kernel of $J$. The right side of above equivalence is separable, so $A_J$ is separable too. Thus we have
\\
\begin {theorem} Local and asymptotic cyclic homology and 
cohomology of $C^*$-algebra $A$ are rigid under deformation quantization by actions of $\Bbb{R}^n$.
\end {theorem}
\begin {proof}
We consider only asymptotic cyclic homology, others are similar. From above discussion, remark 1.5 and theorem 2.1 we have 
\[
HC_*^{\alpha}(A_J)\cong HC_{*+m+n}^{\alpha}(A).
\]
Since $J$ is a skew-symmetric matrix, $m$ is odd if and only $n$ is odd, so always $m+n$ is even.
\end {proof}
\begin {remarks} 
\begin{itemize}
\item[(a)] As an example, we consider n-dimensional noncommutative trous $\Bbb{T}_{\theta}^n$. Even and odd local and asymptotic cyclic homology and cohomology groups of $\Bbb{T}_{\theta}^n$ are isomorphic to  
$\Bbb{C}^{2^{n-1}}$, see 10.2 of [R5] and remark 2.6 (b).
\item[(b)] Theorems 7.5 of [Pu1] and 3.19 of [Pu2] show the inclusion 
$A^{\infty}\rightarrow A$ induces isomorphisms between asymptotic and local cyclic homology and cohomology groups. 
\item[(c)] Let $\alpha$ also denotes the action induced by $\alpha$ on $A_J$, 
subalgebra of smooth elements of this action is $A^{\infty}$ too, theorem 7.1 of [R5]. By part (b), one yields asymptotic and local cyclic (co)homology of smooth algebras are rigid after deformation quantization too.
\end{itemize}
\end {remarks}
\noindent{\bf Acknowledgment.} {I would like to thank Professors R. J. Stanton, H. 
Moscovici and J. Cuntz for their help and guidance. Without their guidance, it could be very difficult for me to complete this work.}

\bibliographystyle {amsalpha}

\begin {thebibliography} {VDN92}
\bibitem [A1]{a} {\sc Abadie, B.} Vector bundles over quantum Heisenberg manifolds. 
{\it Algebraic methods in operator theory.} Brikhauser Boston, Boston, MA, 307-315, 1994.
\bibitem [A2]{a} {\sc Abadie, B.} Generalized fixed-point algebras of certain actions on 
crossed products. {\it Pacific. J. Math}. {\bf 171} (1995), no.1, 1-21.
\bibitem [AEE] {aee} {\sc Abadie, B., Eilers, S., Exel, R.} Morita equivalence for crossed 
products by Hilbert $C^*$-bimodules. {\it Transactions of the 
Amer. Math. Soc}. {\bf 350} (1998), no.8, 3043-3054.
\bibitem [B] {b} {\sc Blackadar, B.} $K$-theory for operator algebras. MSRI Publications, 
Springer-Verlag, New York, 1986.
\bibitem [BGR] {bgr} {\sc Brown, L. G., Green, Ph., Rieffel, M. A.} Stable isomorphism and 
strong Morita equivalence of $C^*$-algebras. {\it Pacific J. Math.} {\bf 71} 
(1977), no. 2, 349-363.
\bibitem [C1] {c} {\sc Cuntz, J.} Generalized homomorphisms between $C^*$-algebras and 
$KK$-theory. {\it Dynamics and processes (Bielefeld, 1981), 31-45, 
Lecture Notes in Math}. {\bf 1031}, {\it Springer, Berlin,} 1983.
\bibitem [C2] {c} {\sc Cuntz, J.} $K$-theory and $C^*$-algebras. {\it Algebraic $K$-theory, 
number theory, geometry and analysis (Bielefeld, 1982), 55-79, Lecture Notes in 
Math}. {\bf 1046}, {\it Springer, Berlin}, 1984.
\bibitem [C3] {c} {\sc Cuntz, J.} A new look at $KK$-theory. {\it $K$-Theory} {\bf 1} 
(1987), no.1, 31-51.
\bibitem [E1] {e} {\sc Exel, R.} A Fredholm operator approach to Morita equivalence. 
{\it $K$-theory} {\bf 7} (1993), no.3, 285-308.
\bibitem [E2] {e} {\sc Exel, R.} Circle actions on $C^*$-algebras, partial automorphisms, 
and a generalized Pimsner-Voiculescu exact sequence. {\it J. Funct. Anal.} {\bf 
122} (1994), no.2, 361-401.
\bibitem [FS] {fs} {\sc Fack, T., Skandalis, G.} Connes' analogue of the Thom 
isomorphism for the Kasparov groups, {\it Inventiones Math.} {\bf 64} (1981), 
7-14.
\bibitem [KhS] {ks} {\sc Khoshkam, M., Skandalis, G.} Toeplitz algebras associated with 
endomorphisms and Pimsner-Voiculescu exact sequences. 
{\it Pacific J. Math.} {\bf 181} (1997), no.2, 315-331.
\bibitem [L] {l} {\sc Loday, J.-L.} Cyclic homology, {\it Grundlehren der math. 
Wissenschaften,} 
Vol. {\bf 310} Springer-Verlag, Brlin, Heidelbelrg, and New York, 1998.
\bibitem [NT1] {nt} {\sc Nest, R., Tsygan, B.} Algebraic index theorem, {\it Comm. 
Math. Phys}.{\bf 172} (1995), 223-262.
\bibitem [NT2] {nt} {\sc Nest, R., Tsygan, B.} Algebraic index theorem for families, 
{\it Adv. Math}. {\bf 113} (1995), 151-205.
\bibitem [P]{p} {\sc Pimsner, M.} A class of $C^*$-algebras generalizing both 
Cuntz-Krieger 
algebras and crossed products by $Z$. {\it Free probability 
theory (Waterloo, ON, 1995), 189-212, Fields Inst. Commun.}
 {\bf 12}, Amer. Math. Soc. Providence, RI, 1997.
\bibitem [PV] {pv} {\sc Pimsner, M., Voiculescu, D.} Exact sequences for $K$-groups and 
Ext-groups of certain cross-product $C^*$-algebras. {\it J. Operator Theory} 
{\bf 4} (1980), no.1, 93-118. 
\bibitem [Pu1]{p}{\sc Puschnigg, M.} Asymptotic cyclic cohomology. {\it Springer 
Lecture Notes in Mathematics}. {\bf 1642} (1996).
\bibitem [Pu2]{p}{\sc Puschnigg, M.} Cyclic homology theories for topological algebras. 
 {\it $K$- theory  Preprint Archives} {\bf 292}.
\bibitem [Pu3]{p}{\sc Puschnigg, M.} Excision in cyclic homology theories. {\it Invent.
 Math}. {\bf 143}, (2001), 249-323.
\bibitem [R1]{r}{\sc Rieffel, M. A.} Morita equivalence for $C^*$-algebras and 
$W^*$-algebras. 
{\it J. Pure Appl. Algebra} {\bf 5} (1974), 51-96.
\bibitem [R2]{r}{\sc Rieffel, M. A.} Induced representations of $C^*$-algebras. 
{\it Advances in Math.} {\bf 13}  (1974), 176-257.
\bibitem [R3]{r}{\sc Rieffel, M. A.} Deformation quantization of Heisenberg manifolds. 
{\it Comm. Math. Phys}. {\bf 122} (1989), no. 4, 531--562.
\bibitem [R4]{r}{\sc Rieffel, M. A.} Noncommutative tori-a case study of noncommutative 
differentiable manifolds. {\it Geometric and topological 
invariants of elliptic operators ( Brunswick, ME, 1988),} 191-211, Contemp. 
Math. Vol. {\bf 105}, Ame. Math. Soc., Providence, RI, 1990.
\bibitem [R5]{r}{\sc Rieffel, M. A.} Deformation quantization for actions of $\Bbb{R}^d$. 
{\it Mem. Amer. Math. Soc}. {\bf 106} (1993), no. 506, x+93 pp.
\bibitem [R6]{r}{\sc Rieffel, M. A.} $K$-groups of $C^*$-algebras deformed by actions 
of $\Bbb{R}^d$. J. {\it Funct. Anal}. {\bf 116} (1993), no.1 199-214.
\bibitem [Ro1]{r}{\sc Rosenberg, J. M.} Algebraic K-theory and its applications. 
{\it Graduate texts in math.} Vol. {\bf 147} Springer-Verlag, Brlin, 
Heidelbelrg, and New York, 1994.
\bibitem [Ro2]{r}{\sc Rosenberg, J. M.} Rigidity of $K$-theory under deformation 
quantization. q-alg/ {\bf 9607021}.
\bibitem [Ro3]{r}{\sc Rosenberg, J. M.} Behavior of $K$-theory under quantization. 
{\it Operator algebras and quantum field theory (Rome, 1996)}, 404-415, 
{\it Internat. Press, Cambridge, MA}, 1997.
\bibitem [S]{s}{\sc Skandalis, G.} Exact sequences for the Kasparov groups of graded 
algebras. 
{\it Canad. J. Math.} {\bf 37} (1985), no.2, 193-216.
\end {thebibliography}

\end {document}